
\input amssym.def
\input amssym
\input epsf
\magnification=1100
\baselineskip = 0.23truein
\lineskiplimit = 0.01truein
\lineskip = 0.01truein
\vsize = 8.5truein
\voffset = 0.2truein
\parskip = 0.10truein
\parindent = 0.3truein
\settabs 12 \columns
\hsize = 5.8truein
\hoffset = 0.2truein

\setbox\strutbox=\hbox{%
\vrule height .708\baselineskip
depth .292\baselineskip
width 0pt}
\font\caps=cmcsc10

\font\bigtenrm=cmr10 at 14pt

\def\sqr#1#2{{\vcenter{\vbox{\hrule height.#2pt
\hbox{\vrule width.#2pt height#1pt \kern#1pt
\vrule width.#2pt}
\hrule height.#2pt}}}}
\def\square{\mathchoice\sqr46\sqr46\sqr{3.1}6\sqr{2.3}4}

\def\G{\Gamma}

\centerline{\bigtenrm SURFACE SUBGROUPS}
\centerline{\bigtenrm OF KLEINIAN GROUPS WITH TORSION}
\tenrm
\vskip 14pt
\centerline{MARC LACKENBY\footnote*{Supported by an EPSRC Advanced Research Fellowship}}
\vskip 18pt
\centerline{\caps 1. Introduction}
\vskip 6pt

A subgroup of a group is known as a {\sl surface subgroup} if it is isomorphic to
the fundamental group of a closed connected orientable surface with positive genus. The {\sl surface subgroup conjecture}
in low-dimensional topology proposes that if the fundamental group of a closed orientable irreducible
3-manifold is infinite, then it contains a surface subgroup. With the recent
solution to the geometrisation conjecture by Perelman ([19], [20], [21]), it suffices to
prove the surface subgroup conjecture for closed orientable hyperbolic 3-manifolds. Gromov
has asked a much more general question [1]: does every infinite, word hyperbolic group
that is not virtually free contain a surface subgroup? In this paper,
we will establish this conjecture for a large class of Kleinian
groups. By a {\sl Kleinian group}, we mean a discrete
subgroup of ${\rm PSL}(2, {\Bbb C})$.

\noindent {\bf Theorem 1.1.} {\sl Any finitely generated, Kleinian group that contains a finite,
non-cyclic subgroup either is finite or virtually free or contains a surface subgroup.}

The orientation-preserving isometry group of hyperbolic 3-space ${\Bbb H}^3$ is ${\rm PSL}(2, {\Bbb C})$
and the quotient of ${\Bbb H}^3$ by a Kleinian group $\G$ is an orientable hyperbolic
3-orbifold $O$. The hypothesis that $\G$ contains a finite,
non-cyclic subgroup is equivalent to the statement that the
singular locus of $O$ contains a vertex. 
When the orbifold $O$ is non-compact, Theorem 1.1 has been known for many years.
(We will recall the proof later in the Introduction.)
Thus, the heart of Theorem 1.1 is a statement about closed hyperbolic
3-orbifolds.

Whether or not a group contains a surface subgroup depends only on
its commensurability class. In other words, if two groups contain
finite index subgroups that are isomorphic, then either they both
contain a surface subgroup or neither does. It is a theorem of
the author, Long and Reid [14] that any arithmetic Kleinian group
is commensurable with one that contains ${\Bbb Z}/2 \times {\Bbb Z}/2$.
Moreover, arithmetic Kleinian groups have finite co-volume and so are
neither finite nor virtually free.
Thus, we obtain the following corollary to Theorem 1.1.

\noindent {\bf Theorem 1.2.} {\sl Any arithmetic Kleinian group
contains a surface subgroup.}

\vfill\eject
Some authors restrict the term `surface subgroup' to mean the fundamental
group of a closed connected orientable surface with genus at least two. We have chosen
not to do so here, in order to state the results most succinctly. However,
in the main case we consider, when the Kleinian group $\G$ is co-compact,
the surface subgroup cannot be ${\Bbb Z}\times {\Bbb Z}$. It should also
be pointed out that Cooper, Long and Reid [7] proved that, when $\G \backslash{\Bbb H}^3$
has finite volume but is not compact, then $\G$ contains the fundamental
group of a closed connected orientable surface with genus at least two.

The surface subgroup conjecture is the first in a hierarchy of
increasingly strong conjectures. For convenience, we state
these only for a closed orientable hyperbolic \break 3-manifold $M$. The virtually
Haken conjecture asserts that $M$ has a finite cover that is Haken, that is, it contains
a closed embedded orientable incompressible surface with positive genus. 
By the loop theorem, such a surface is $\pi_1$-injective,
and so gives rise to a surface subgroup of $\pi_1(M)$. The positive
virtual $b_1$ conjecture [25] proposes that $M$ has a finite cover with
positive first Betti number. The infinite virtual $b_1$ conjecture
asserts that $M$ has finite covers $\tilde M$ where $b_1(\tilde M)$
is arbitrarily big. The largeness conjecture proposes that $\pi_1(M)$
is large, that is, it contains a finite index subgroup that admits
a free non-abelian quotient. In another direction, the virtual fibering conjecture [25]
asserts that $M$ has a finite cover that fibres over the circle.
The proof of Theorem 1.1 falls into
a general programme initiated by the author [11] for attacking some of these
conjectures, by relating the problem to an analysis of Cheeger
constants.

Recall that the {\sl Cheeger constant} of a closed Riemannian manifold
$M$ is defined to be
$$h(M) = \inf_S \left \{ {{\rm Area}(S) \over \min \{ {\rm Volume}(M_1), {\rm Volume}(M_2) \} } \right \},$$
where $S$ ranges over all smooth codimension one submanifolds that separate
$M$ into submanifolds $M_1$ and $M_2$. The {\sl Cheeger constant} of an
open Riemannian manifold $M$ with infinite volume is
$$h(M) = \inf_S \left \{ {{\rm Area}(S) \over {\rm Volume}(M_1) } \right \},$$
where $S$ ranges over all smooth codimension one
submanifolds that bound a compact submanifold $M_1$.

A central conjecture, due to Lubotzky and Sarnak [17], proposes that
any closed orientable hyperbolic 3-manifold has a sequence of finite
covers $M_i$ such that $h(M_i) \rightarrow 0$. A theorem of the author,
Long and Reid [14] gives evidence for the reach of this conjecture.
It states that if the Lubotzky-Sarnak conjecture is true for every
closed hyperbolic 3-manifold, then any co-compact Kleinian group containing
${\Bbb Z}/2 \times {\Bbb Z}/2$ is large. In addition, a theorem of
the author [11] states that the Lubotzky-Sarnak conjecture, plus another
conjecture called the Heegaard gradient conjecture, together imply
the virtually Haken conjecture. Unfortunately, the Lubotzky-Sarnak
conjecture remains wide open at present, despite a considerable amount
of positive evidence [8]. However, the following theorem of the author, Long
and Reid [15] establishes a version of the conjecture for infinite
covers.

\noindent {\bf Theorem 1.3.} {\sl Any closed orientable hyperbolic
3-manifold has a sequence of infinite covers $M_i$ such that $h(M_i) \rightarrow 0$.}

This is a consequence of the following theorem of Bowen [2]. 
In Section 5, we give an alternative interpretation of part of Bowen's proof.

\noindent {\bf Theorem 1.4.} {\sl Let $\G$ be a co-compact lattice in ${\rm PSL}(2, {\Bbb C})$.
Then $\G$ contains a sequence of finitely generated, free, convex co-compact
subgroups $\G_i$ such that the Hausdorff dimension of the limit set
of $\G_i$ tends to $2$.}

Briefly, the argument in [15] proving Theorem 1.3 from Theorem 1.4 runs as follows. Let $M_i = \G_i \backslash {\Bbb H}^3$
and let $D_i$ be the Hausdorff dimension of the limit set of $\G_i$. A theorem of Sullivan [24] gives
that $\lambda_0(M_i) = D_i (2 - D_i)$, provided $D_i \geq 1$, where
$\lambda_0(M_i)$ is the infimum of the spectrum of the Laplacian on $L_2$ functions
on $M_i$. Then Cheeger's inequality [4] gives that $h(M_i) \leq \lambda_0(M_i)^2/4 \rightarrow 0$.

We now explain why having Cheeger constants tend to zero is so useful,
and place the proofs of Theorems 1.1 and 1.2 into context.
As mentioned above, the author proposed an approach in [11] to the virtually Haken conjecture.
If one can find a nested sequence of finite regular covers $M_i$ of a closed orientable
irreducible 3-manifold such that
$h(M_i) \rightarrow 0$ but where the Heegaard genus of $M_i$ grows linearly
as a function of the covering degree, then Theorem 1.7 of [11] states that
the manifolds $M_i$ are eventually Haken. (The Heegaard genus of a closed orientable 3-manifold
is the minimal integer $g$ such that the manifold can be expressed as the union of
two genus $g$ handlebodies glued homeomorphically along their boundaries.)
In [13], the author applied these ideas to
closed orientable hyperbolic 3-orbifolds $O$ with non-empty singular locus. 
There it was shown that $O$ has a sequence of finite covers $O_i$
where the number of simple closed curve components of the singular locus grows linearly as a function of
the covering degree. If the singularity order of these components is a
multiple of a prime $p$, then it was also shown that the rank of $H_1(O_i; {\Bbb Z}/p)$
is at least the number of such components. Hence,
the rank of $H_1(O_i; {\Bbb Z}/p)$ grows linearly as a function of the
covering degree. But the Heegaard genus is at least the rank of
$H_1(O_i; {\Bbb Z}/p)$. So, if, in addition, we can arrange that
$h(O_i) \rightarrow 0$, then Theorem 1.7 in [11] implies that $O$
is virtually Haken. In fact, Theorem 1.7 in [13] gives much more:
under these hypotheses, $\pi_1(O)$ is large. This was proved using Theorem 1.1
of [12], which shows how largeness follows from linear growth of the rank of
mod $p$ homology, together with Cheeger constants tending to zero. This demonstrates why
3-orbifolds with non-empty singular locus appear to be much more tractable
than 3-manifolds in questions related to the virtually Haken conjecture.
In [14], the author, Long and Reid developed the technology in [13] and
applied it to the case where the singular locus of $O$ has a vertex.
We particularly focused on the case where every singular vertex of $O$ 
has local group ${\Bbb Z}/2 \times {\Bbb Z}/2$ and there is at least one such vertex.
It was shown that if $|O_i|$ is any sequence of finite covers of the
underlying manifold $|O|$ and $O_i$ is the corresponding cover of $O$, then $H_1(O_i; {\Bbb Z}/2)$ grows linearly
as a function of the covering degree. If, in addition, $h(|O_i|) \rightarrow 0$,
then $\pi_1(O)$ is large.  In [14], the author, Long and Reid also showed that any
arithmetic Kleinian group is commensurable with one containing 
${\Bbb Z}/2 \times {\Bbb Z}/2$. Thus, these methods apply to arithmetic
hyperbolic 3-manifolds.

In the arithmetic case, the largeness conjecture (and hence the surface subgroup
conjecture) is already known
to be true in certain cases, which we now describe. By definition, any arithmetic Kleinian group $\G$
is commensurable with the units in an order of a quaternion algebra $D$
defined over a number field $L$ with a unique complex embedding and such that
$D$ ramifies for all real embeddings of $L$. When $L$ has a subfield of
index 2, it was shown by Labesse and Schwermer [10] that $\G$ has a finite index
subgroup with positive first Betti number (see also [17]). The same conclusion was also
established by Clozel [5] under the hypothesis that, for each prime $\nu$ of
$L$ for which $D$ ramifies at $L_\nu$, $L_\nu$ contains no quadratic extension
of ${\Bbb Q}_p$, where $p$ is the rational prime that $\nu$ divides. In all the
cases where the arithmetic methods of [10] and [5] could produce a
finite index subgroup with positive first Betti number, the author, Long and
Reid in [14] were able to show that $\G$ is in fact large.

The proof of Theorem 1.1 proceeds as follows. Let $\G$ be the finitely generated Kleinian group that contains a finite,
non-cyclic subgroup. Let $O$ be the 3-orbifold $\G \backslash {\Bbb H}^3$. If $O$ is
not closed, then Theorem 1.1 is well known. The argument runs as follows. Selberg's lemma implies
that $O$ has a finite cover that is an orientable irreducible 3-manifold $M$. This has
finitely generated fundamental group, and so Scott's theorem [22]
gives that $\pi_1(M)$ is isomorphic to the fundamental group of
a compact orientable irreducible 3-manifold $M'$ with non-empty boundary. If the boundary of
$M'$ is completely compressible, then $M'$ is a 3-ball or handlebody and so
$\pi_1(M)$ is trivial or free. If not, then $M'$ contains a closed orientable properly
embedded incompressible surface with positive genus. This surface is
$\pi_1$-injective (see Corollary 6.2 in [9] for example). Hence,
$\pi_1(M)$ contains a surface subgroup. 

Therefore, let  us assume that $O$ is a closed 3-orbifold. We first show that $O$ has a finite cover
$\tilde O$ such that the underlying manifold $|\tilde O|$ of $\tilde O$ has infinite
fundamental group, and where $\tilde O$ contains a singular vertex. (We also
need some other conditions on the singular locus, which we omit here.) In the
case where $\pi_1(O)$ contains ${\Bbb Z}/2 \times {\Bbb Z}/2$, this was
proved by the author, Long and Reid in [14]. But if $\pi_1(O)$ does not contain
${\Bbb Z}/2 \times {\Bbb Z}/2$, then an alternative proof is required,
which is given in Section 4. Now, it
is a well known consequence of Perelman's solution to the geometrisation
conjecture that if $M$ is a closed orientable 3-manifold with infinite fundamental
group, then either $M$ has a finite cover with
positive first Betti number or $M$ is hyperbolic. We apply this to $M = |\tilde O|$.
If $M$ has a finite cover with positive first Betti number, then $\pi_1(\tilde O)$ clearly contains a surface
subgroup. Thus, we may assume that $|\tilde O|$ is hyperbolic. We then use
Theorem 1.3 to give a sequence of covers $|O_i|$ of $|\tilde O|$
such that $h(|O_i|) \rightarrow 0$. Let $O_i$ be the corresponding
covers of $\tilde O$. These are infinite covers, and so the 
technology developed by the author in [12] to prove largeness does not directly apply. Nevertheless,
we can adapt it to prove the existence of a surface subgroup. 

It is intriguing
that the proofs of Theorems 1.1 and 1.2, which are purely statements about hyperbolic
3-orbifolds, should require the solution to the geometrisation conjecture in such an
essential way.

\noindent {\sl Acknowledgement.} The author would like to thank the referees of this
paper for their helpful suggestions, which have improved it considerably.

\vskip 18pt
\centerline{\caps 2. Preliminaries on orbifolds}
\vskip 6pt

In this section, we briefly introduce some standard terminology and results
about orbifolds. We assume some familiarity with their basic theory. 
In particular, we will not define covering spaces and the
fundamental group of orbifolds. (See, for example [23] or [6].)
We let $\pi_1(O)$ denote the orbifold fundamental group of an orbifold $O$,
and let $|O|$ denote the underlying topological space of $O$.

Each point $x$ of an orientable $n$-dimensional orbifold $O$ (without boundary) has an open neighbourhood 
modelled on $L(x) \backslash {\Bbb R}^n$, for some finite subgroup $L(x)$ of ${\rm SO}(n)$, where
$x$ is identified with the image of the origin. Similarly, when $O$ is an orientable $n$-dimensional
orbifold, possibly with boundary, each point $x$ of $O$ has a local model of the form $L(x) \backslash
{\Bbb R}^n$ or $L(x) \backslash {\Bbb R}^n_+$, where ${\Bbb R}^n_+ = \{ (x_1, \dots, x_n) \in {\Bbb R}^n : x_n \geq 0 \}$.
Again, $L(x)$ is a finite subgroup of ${\rm SO}(n)$, and in the case where it acts on
${\Bbb R}^n_+$, it is required to preserve the $x_n$ co-ordinate. The group $L(x)$
is known as the {\sl local group} of $x$. The set of points
$x$ where $L(x)$ is non-trivial is known as the {\sl singular locus}
of $O$. We denote this set by ${\rm sing}(O)$. 

One can compute the fundamental group $\pi_1(O)$ as follows.
The points $x$ where $L(x)$ is non-trivial and fixes some codimension 2 subspace of ${\Bbb R}^n$ or ${\Bbb R}^n_+$
form the {\sl codimension 2 singular locus}. The local group $L(x)$ of any such $x$
is cyclic. Let $C_1, \dots, C_r$ denote the
components of the codimension 2 singular locus. If two points lie in
the same component, they have isomorphic local groups. Pick a basepoint
$b$ in $O - {\rm sing}(O)$. For each component $C_i$,
pick a path with interior in $O - {\rm sing}(O)$ running from $b$ to some point $x_i$ in $C_i$.
(This is possible because ${\rm sing}(O)$ does not separate $O$.)
Let $\mu_i$ be the element of $\pi_1(O - {\rm sing}(O))$ represented by the loop
which runs from $b$ along this path almost as far as $x_i$, encircles $C_i$
at $x_i$, and then returns to $b$ along the path. Let $n_i$ denote the
order of the local group at $x_i$. Then
the fundamental group $\pi_1(O)$ is isomorphic to
$$\pi_1(O - {\rm sing}(O),b) / 
\langle \! \langle \mu_1^{n_1}, \dots, \mu_r^{n_r} \rangle \! \rangle.$$
We term this a {\sl meridional presentation} for $\pi_1(O)$.

We now consider a compact orientable 3-orbifold $O$, possibly with boundary. The finite subgroups of
${\rm SO}(3)$ are either trivial, cyclic, dihedral, or
isomorphic to $A_4$, $S_4$ or $A_5$. (We include ${\Bbb Z}/2 \times {\Bbb Z}/2$ as 
a dihedral group.) The quotient of ${\Bbb R}^3$ (respectively, ${\Bbb R}^3_+$) by any such subgroup
is again homeomorphic to ${\Bbb R}^3$ (respectively, ${\Bbb R}^3_+$). Thus, the underlying space $|O|$ of $O$ is a 3-manifold.
This observation will be particularly important in this paper, since we will apply
extensive machinery from 3-manifold theory to $|O|$.

The points in ${\rm sing}(O)$ where the local group is non-trivial and cyclic
form the codimension 2 singular locus, which is a 1-manifold. In any given component 
of this 1-manifold, the points all have the same local group up to isomorphism. Its order is the {\sl
singularity order} of this component.

The points in ${\rm sing}(O)$ where the local group is non-trivial
and non-cyclic form trivalent isolated vertices in ${\rm sing}(O)$. We will mostly
be concerned with vertices having dihedral local group of order $2n$.
Emanating from these vertices are two edges of ${\rm sing}(O)$ with order $2$
and an edge with order $n$.

For any prime $p$, we let ${\rm sing}_p(O)$ denote the union of the 
components of the codimension 2 singular locus with singularity order that is a multiple of $p$.
Thus, ${\rm sing}_p(O)$ is a collection of simple closed curves and
intervals. Let $\overline{{\rm sing}_p(O)}$ denote its closure.

\noindent {\bf Lemma 2.1.} {\sl Let $O$ be a closed orientable 3-orbifold
with singular locus a (possibly empty) collection of simple closed curves.
Then $\pi_1(O)$ has a finite presentation $\langle X | R \rangle$ such that
$|R| - |X|$ is at most the number of components of ${\rm sing}(O)$.}

\noindent {\sl Proof.} This is contained in the proof of Theorem 4.1 in [13] and in the proof of
Lemma 8.3 in [14], but we repeat the argument here. The manifold $O - {\rm sing}(O)$ has the same
fundamental group as $O - {\rm int}(N({\rm sing}(O))$ which is a compact
orientable 3-manifold with boundary a (possibly empty) collection of tori.
Hence, $\pi_1(O - {\rm sing}(O))$ has a presentation with the same number
of generators as relations. To obtain $\pi_1(O)$, we use a meridional presentation.
Starting from $\pi_1(O - {\rm sing}(O))$, we add a relation for
each component of ${\rm sing}(O)$, and the lemma follows
immediately. $\square$

For any group, space or orbifold $X$ and any prime $p$, we let $b_1(X)$
be its first Betti number and let $d_p(X)$ be the dimension of $H_1(X; {\Bbb Z}/p)$
as a vector space over ${\Bbb Z}/p$. The above formula for $\pi_1(O)$
also has the following consequence for homology (which also appears as
Proposition 3.1 in [14]).

\noindent {\bf Lemma 2.2.} {\sl Let $O$ be a compact orientable 3-orbifold,
possibly with boundary, and let $p$ be a prime. Then $d_p(O) \geq b_1(\overline{{\rm sing}_p(O)})$.}

\noindent {\sl Proof.} By the above formula for $\pi_1(O)$, we deduce that 
$$H_1(O; {\Bbb Z}/p) \cong H_1(O - {\rm sing}(O); {\Bbb Z}/p) /
\langle \! \langle \mu_1^{n_1}, \dots, \mu_r^{n_r} \rangle \! \rangle.$$
When $n_i$ is coprime to $p$, quotienting $H_1(O - {\rm sing}(O); {\Bbb Z}/p)$ by
$\mu_i^{n_i}$ has the same effect as quotienting by $\mu_i$.
On the other hand, when $n_i$ is a multiple of $p$, then quotienting
by $\mu_i^{n_i}$ has no effect. Thus, if we let $M$ denote the
3-manifold obtained by removing an open regular neighbourhood of
$\overline{{\rm sing}_p(O)}$, then $d_p(O) = d_p(M)$.
Now, it is a well known consequence of Poincar\'e duality
that $d_p(M) \geq {1 \over 2} d_p(\partial M)$. But $d_p(\partial M)$
is clearly at least $2 b_1(\overline{{\rm sing}_p(O)})$, and
the lemma follows. $\square$

Let us now start the proof of Theorem 1.1. The given finite non-cyclic subgroup of 
the finitely generated Kleinian group $\G$ is either 
dihedral or isomorphic to $A_4$, $S_4$ or $A_5$. In fact, $A_4$, $S_4$ and $A_5$
contain dihedral subgroups. So, we may assume that the subgroup is dihedral
of order $2n$. If $n$ is even, then the subgroup contains ${\Bbb Z}/2 \times {\Bbb Z}/2$.
If $n$ is odd, then the subgroup contains a dihedral subgroup of order
$2p$, where $p$ is an odd prime. We may therefore suppose that
$\G$ contains ${\Bbb Z}/2 \times {\Bbb Z}/2$ or a dihedral group of order
$2p$, where $p$ is an odd prime. The proof divides according to these
cases, which we deal with in Sections 3 and 4 respectively. 

\vskip 18pt
\centerline{\caps 3. The ${\Bbb Z}/2 \times {\Bbb Z}/2$ case}
\vskip 6pt

In this section, we suppose that the finitely generated Kleinian group $\G$ contains a subgroup isomorphic
to ${\Bbb Z}/2 \times {\Bbb Z}/2$. As mentioned in the Introduction, we may also
assume that the quotient $\G \backslash {\Bbb H}^3$ is a closed orientable
3-orbifold $O$.

The following is a slight strengthening of Theorem 8.1 in [14], with essentially
the same proof.

\noindent {\bf Theorem 3.1.} {\sl Let $O$ be a closed orientable hyperbolic 3-orbifold
such that $\pi_1(O)$ contains ${\Bbb Z}/2 \times {\Bbb Z}/2$. Then 
$O$ has a finite cover $\tilde O$ such that
\item{1.} every arc and simple closed curve of ${\rm sing}(\tilde O)$ has order 2;
\item{2.} every vertex of ${\rm sing}(\tilde O)$ has local group
${\Bbb Z}/2 \times {\Bbb Z}/2$ and there is at least one such vertex;
\item{3.} $\pi_1(|\tilde O|)$ is infinite.}

The final conclusion is the difficult part. It was proved using the Golod-Shafarevich
inequality. Recall from [16] that this asserts that if $G$ is a group with
the finite presentation $\langle X|R \rangle$ such that
$d_p(G)^2/4 - d_p(G) + |X| - |R| > 0$, then $G$ is infinite.

A well-known consequence of the geometrisation conjecture 
(see the Appendix in [14] for example) is that if $M$ is a closed
orientable 3-manifold with infinite fundamental group then either $M$ has a finite cover with
positive first Betti number or $M$ is hyperbolic. Let us
suppose first that $|\tilde O|$ has a finite cover $|\tilde O'|$ with positive first
Betti number. Then, $\tilde O'$ is a closed hyperbolic 3-orbifold
with positive first Betti number. Its fundamental group therefore
contains a surface subgroup. In fact, much more can be said in this
case, via Theorem 9.5 in [14], which is as follows.

\noindent {\bf Theorem 3.2.} {\sl Let $\tilde O$ be a closed
orientable hyperbolic 3-orbifold such that $\pi_1(\tilde O)$ contains ${\Bbb Z}/2 \times {\Bbb Z}/2$. 
Suppose that $|\tilde O|$ has a finite
cover with positive first Betti number. Then $\pi_1(\tilde O)$ is
large.}

Thus, we may assume that $|\tilde O|$ is hyperbolic. By Theorem 1.3,
$|\tilde O|$ has a sequence of infinite covers $|O_i|$ where $h(|O_i|) \rightarrow 0$.
Let $O_i$ be the corresponding covers of $\tilde O$.

Let us fix a 1-vertex triangulation $T$ of $|\tilde O|$. For convenience,
we may suppose that its vertex is one of the vertices of ${\rm sing}(\tilde O)$ 
with local group ${\Bbb Z}/2 \times {\Bbb Z}/2$.
The inverse image of this triangulation in $|O_i|$ is again a
triangulation $T_i$. Let $X_i$ be its 1-skeleton. 

It is convenient to work with the Cheeger constant of $X_i$. Recall
that this is defined as follows. For a graph $X$,
let $V(X)$ and $E(X)$ denote its vertices and edges. For $A \subset V(X)$,
let $\partial A$ denote the set of edges with one endpoint in $A$
and one endpoint not in $A$. Define the {\sl Cheeger constant} of a
graph $X$ where $V(X)$ is infinite to be
$$h(X) = \inf \left \{ {|\partial A| \over |A| } : \emptyset \not= A \subset V(X) \hbox{ and } A
\hbox{ is finite} \right \}.$$
It is a standard fact, first observed by Brooks (see the proof of Theorem 1 in [3] for example),
that there is a constant $k \geq 1$ such that
$$k^{-1} \ h(|O_i|) \leq h(X_i) \leq k \ h(|O_i|).$$
So, $h(X_i) \rightarrow 0$. For each $i$, pick a non-empty finite set $A_i$ of vertices in $X_i$
such that $|\partial A_i|/|A_i| \rightarrow 0$. We now use
this set to create a closed orientable 2-dimensional sub-orbifold $S_i$ in $O_i$ bounding a
compact 3-dimensional sub-orbifold $N_i$.

In the interior of each 1-cell of $T$, pick a point,
which we declare to be its midpoint. In each triangular face of $T$,
pick three properly embedded arcs, each joining two of the three midpoints
on its boundary. We may ensure that these arcs have disjoint interiors. We term these {\sl specified normal arcs}.
We may suppose that these arcs are disjoint from the singular locus of $\tilde O$.
For each tetrahedron in $T$, consider all the simple closed curves
in its boundary which are the union of specified normal arcs, and
which are transverse to the 1-skeleton of the tetrahedron.
There are seven such curves, consisting
of four triangles and three quadrilaterals. For each
such curve, pick a properly embedded disc that it bounds
in the tetrahedron. We term this a {\sl specified normal disc}.
We may suppose that these discs intersect
the singular locus of $\tilde O$ transversely, and are disjoint
from the vertices of ${\rm sing}(\tilde O)$. The specified normal arcs and
discs lift to arcs and discs in $|O_i|$, which we also call
{\sl specified normal arcs and discs}. We will use these to construct $S_i$.

Consider any face of $T_i$ that contains an edge of $\partial A_i$.
Then it contains precisely two edges of $\partial A_i$. Join these
by a specified normal arc. Now consider a tetrahedron in $T_i$
that contains an edge of $\partial A_i$. Either three or four edges
of this tetrahedron belong to $\partial A_i$, and these have
been joined by specified normal arcs, yielding a triangle or quadrilateral curve.
Fill this in with a specified normal disc. The union of these
discs is the surface $|S_i|$. Because the normal discs are transverse to
${\rm sing}(\tilde O)$ and disjoint from the vertices of ${\rm sing}(\tilde O)$,
$|S_i|$ is the underlying space of an orbifold $S_i$.
Note that  $S_i$ separates $O_i$ into two sub-orbifolds.
One of these contains $A_i$ and is compact. We denote
it by $N_i$. Note also that the expression of $|S_i|$ as a union
of normal triangles and quadrilaterals determines a cell structure for it.

\noindent {\sl Claim 1.} There are constants $k_0$, $k_1$, $k_2 > 0$,
independent of $i$, such that the number of 0-cells, 1-cells and 2-cells of $|S_i|$
is at most $k_0 |\partial A_i|$, $k_1 |\partial A_i|$ and
$k_2 |\partial A_i|$ respectively.

The 0-cells of $|S_i|$ are precisely the midpoints of the edges in
$\partial A_i$. So, we may set $k_0 = 1$. The number of 1-cells of
$|S_i|$ is at most the number of
0-cells of $|S_i|$ times half the maximal valence of
any 0-cell of $|S_i|$. But this valence is bounded above, independent of $i$, since
it is at most the maximal number of faces of $T$
running over any 1-cell of $T$ (counted with multiplicity). Finally, the number
of 2-cells of $|S_i|$ is at most $(2/3)$ times the number
of 1-cells, since each 2-cell is a triangle or quadrilateral.
This proves the claim.

\noindent {\sl Claim 2.} There is a constant $k_3 > 0$, independent of
$i$, such that $d_2(S_i) \leq k_3 |\partial A_i|$.

Let $k_4$ be the maximal number of intersection points between
the singular locus of $\tilde O$ and any specified normal disc in $T$.
For any component $S'_i$ of $S_i$, $\pi_1(S'_i)$ is obtained from
$\pi_1(S'_i - {\rm sing}(S'_i))$ by quotienting by the normal subgroup
generated by squares of the meridians of ${\rm sing}(S'_i)$.
Hence, $d_2(S_i) = d_2(S_i - {\rm sing}(S_i))$. But 
$d_2(S_i - {\rm sing}(S_i))\leq d_2(|S_i|) + |S_i \cap {\rm sing}(O_i)|$,
since $H_1(|S_i|; {\Bbb Z}/2)$ is obtained from $H_1(S_i - {\rm sing}(S_i); {\Bbb Z}/2)$
by quotienting out the subgroup generated by the meridians of the singular locus.
Now, $d_2(|S_i|)$ is at most the number of 1-cells of $|S_i|$.
And $|S_i \cap {\rm sing}(O_i)|$ is at most $k_4$ times the
number of 2-cells of $|S_i|$. Thus, the claim follows from Claim 1.

\noindent {\sl Claim 3.} There is a constant $k_5 > 0$, independent of $i$,
such that $d_2(N_i) \geq |A_i|/2 - k_5 |\partial A_i|$.  

By Lemma 2.2, $d_2(N_i) \geq b_1(\overline{{\rm sing}_2(N_i)})$. Now, $\overline{{\rm sing}_2(N_i)}$
is a collection of simple closed curves and graphs. Let $G_i$
be the union of the graph components. The vertices of $G_i$ come in two types. 
Those in the interior of $N_i$ have valence 3, and those on $S_i$ have valence 1.
The number of valence 1 vertices is $|S_i \cap G_i|$,
which is at most $k_4 k_2 |\partial A_i|$. Now,
$b_1(G_i) > -\chi(G_i)$, which is the sum, over each 
vertex $v$ of $G_i$, of $({\rm val}(v)/2) -1$, where ${\rm val}(v)$
is the valence of $v$. The trivalent vertices of $G_i$ each
contribute $1/2$ to this sum, and we have arranged that there is
one of these vertices at each point of $A_i$. So,
$$d_2(N_i) \geq b_1(\overline{{\rm sing}_2(N_i)}) \geq b_1(G_i) \geq |A_i|/2 - (k_4 k_2/2) |\partial A_i|.$$
Setting $k_5 = k_4 k_2/2$ proves the claim.

It is here, in Claim 3, that we are making crucial use of the vertices of the
singular locus. This is why the main theorem requires the existence of
a finite subgroup that is non-cyclic.

By Claims 2 and 3, when $|\partial A_i| / |A_i|$ is sufficiently
small, $d_2(N_i) > d_2(S_i)$. Thus, the kernel of 
$H^1(N_i; {\Bbb Z}/2) \rightarrow H^1(S_i; {\Bbb Z}/2)$
is non-trivial. Let $\alpha_1$ be a non-trivial element of
this kernel. This is the image of some element $\alpha_2 \in H^1(N_i, S_i; {\Bbb Z}/2)$.
By excision, this is isomorphic to $H^1(O_i, O_i - N_i; {\Bbb Z}/2)$.
Let $\alpha_3 \in H^1(O_i, O_i - N_i; {\Bbb Z}/2)$ be the element
corresponding to $\alpha_2$. This is sent to a non-trivial
element $\alpha_4 \in H^1(O_i; {\Bbb Z}/2)$. The image of
$\alpha_4$ in $H^1(O_i - N_i; {\Bbb Z}/2)$ is trivial.
Let $\tilde O_i$ be the 2-fold cover of $O_i$ dual to
$\alpha_4$. Since the image of $\alpha_4$ in $H^1(O_i - N_i; {\Bbb Z}/2)$
is trivial, $\tilde O_i$ contains two copies
of $O_i - N_i$. In particular, $|\tilde O_i|$ has at least two ends. 

\vfill\eject
\centerline{
\epsfxsize=2.5in
\epsfbox{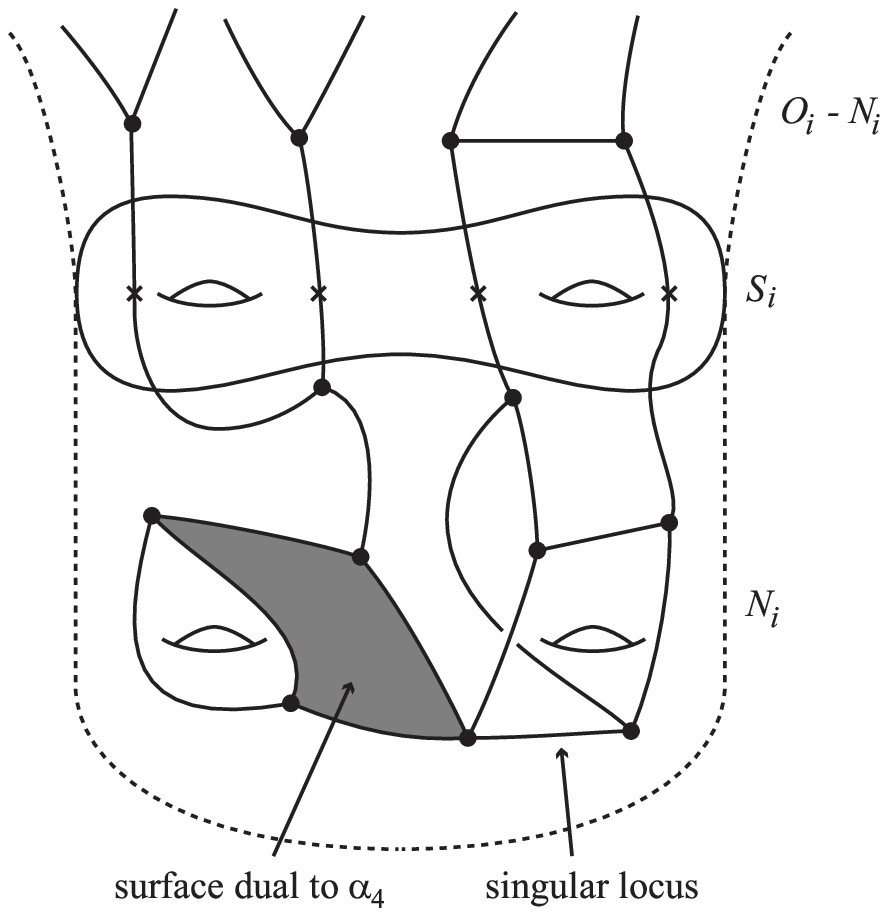}
}
\vskip 6pt
\centerline{Figure 1.}

Selberg's lemma implies that $O$ is finitely covered by a hyperbolic 3-manifold $M$.
Let $M_i$ be the cover of $O$ corresponding to $\pi_1(\tilde O_i) \cap \pi_1(M)$.
This is a finite cover of $\tilde O_i$. Hence, it also has at least two ends. It is irreducible because it is
a hyperbolic 3-manifold. 
We now use the following simple lemma, applied to $M_i$, to complete the proof
of Theorem 1.1, in the case where $\G$ contains ${\Bbb Z}/2 \times
{\Bbb Z}/2$.

\noindent {\bf Lemma 3.3.} {\sl Let $M$ be an orientable irreducible open 3-manifold with
at least two ends. Then $\pi_1(M)$ contains a surface subgroup.}

\noindent {\sl Proof.} Let $S$ be a closed orientable surface separating two ends of $M$.
This surface may admit a {\sl compression disc}, which is a disc $D$ embedded in $M$
such that $D \cap S = \partial D$ and such that $\partial D$ does not bound a disc in $S$.
If so, we {\sl compress} $S$, which is the operation of removing a regular neighbourhood
of $\partial D$ in $S$, and attaching two parallel copies of $D$.
Compress $S$ as much as possible to a surface $\overline S$. Some component of $\overline S$ still
separates two ends of $M$. Hence, it cannot be a 2-sphere,
since $M$ is irreducible. By the loop theorem,
$\overline S$ is $\pi_1$-injective (see Corollary 6.2 in [9]).
So, $\pi_1(M)$ contains a surface subgroup. $\square$

\vskip 18pt
\centerline{\caps 4. The dihedral case}
\vskip 6pt

In this section, we suppose that the finitely generated Kleinian group $\G$ contains a dihedral
group $D$ of order $2p$, where $p$ is an odd prime. We start with
the following theorem, which is an analogue of Theorem 3.1.

\vfill\eject
\noindent {\bf Theorem 4.1.} {\sl Let $O$ be a closed orientable
hyperbolic 3-orbifold such that $\pi_1(O)$ contains a dihedral group
$D$ of order $2p$, where $p$ is an odd prime. Then $O$ has a
finite cover $\tilde O$ such that 
\item{1.} every arc and simple closed curve of ${\rm sing}(\tilde O)$
has order $2$ or $p$;
\item{2.} the vertices of ${\rm sing}(\tilde O)$ each have local group isomorphic to $D$,
and there is at least one such vertex;
\item{3.} $\pi_1(|\tilde O|)$ is infinite.

\noindent In addition, there is a regular finite cover $\tilde O' \rightarrow \tilde O$,
with degree a power of $2$, such that ${\rm sing}(\tilde O')$ is a non-empty collection of simple closed
curves of order $p$.}

\noindent {\sl Proof.} By Selberg's lemma, $O$ admits a finite regular
manifold cover $M$. Let $O_1$ be the covering space of $O$ corresponding
to the subgroup $\pi_1(M)D$ of $\pi_1(O)$. Since $\pi_1(O_1)$ contains $D$, ${\rm sing}(O_1)$
has a vertex with local group $D$. Since $O_1$ is the quotient
of $M$ by a group of covering transformations isomorphic to $D$, every
arc and circle of ${\rm sing}(O_1)$ has order $2$ or $p$. Moreover,
every vertex of ${\rm sing}(O_1)$ has local group isomorphic to $D$.

Let $m = 10$. Theorem 4.6 of [14] states that there is a surjective homomorphism
$\phi \colon \pi_1(O_1) \rightarrow {\rm PSL}(2, p_1) \times \dots \times
{\rm PSL}(2, p_m) = G$, where $p_1, \dots, p_m$ is a collection of primes.
Moreover, if $\pi_i \colon G \rightarrow {\rm PSL}(2, p_i)$ is
projection onto the $i$th factor, then the kernel of $\pi_i \phi$ is
torsion free, for each $i$.

Let $C$ be the index 2 subgroup of $D$ isomorphic
to ${\Bbb Z}/p$. Let $C_i$ (respectively, $D_i$) be the image of
$C$ (respectively, $D$) in ${\rm PSL}(2, p_i)$ under $\pi_i \phi$. Since
the kernel of $\pi_i \phi$ is torsion free, $C_i$
(respectively, $D_i$) is a copy of $C$ (respectively, $D$).

Let $O_2$ be the covering space of $O_1$ corresponding to
$\phi^{-1}(D_1 \times \dots \times D_m)$. Note that $\pi_1(O_2)$ contains
$D$. Hence, ${\rm sing}(O_2)$ contains a vertex with local group $D$.
Since $O_2$ finitely covers $O_1$, every arc and circle of ${\rm sing}(O_2)$
has order $2$ or $p$. Also, every vertex of ${\rm sing}(O_2)$ has
local group isomorphic to $D$.

Now, $\phi^{-1}(C_1 \times \dots \times C_m)$ is a normal subgroup of 
$\phi^{-1}(D_1 \times \dots \times D_m) = \pi_1(O_2)$ with quotient $({\Bbb Z}/2)^m$.
Hence, $d_2(O_2) \geq m = 10$.

Emanating from each vertex of ${\rm sing}(O_2)$, there are two edges with order $2$
and one edge with order $p$. So, $\overline{{\rm sing}_2(O_2)}$ is a collection
of simple closed curves. Pick a component $V$ of $\overline{{\rm sing}_2(O_2)}$ that
contains a vertex. Let $\mu$ be one of its meridians.
Let $\langle \! \langle \mu \rangle \! \rangle$ be the
subgroup of $\pi_1(O_2)$ normally generated by $\mu$.
Let $\pi_1(O_2)^2$ denote the subgroup of $\pi_1(O_2)$ generated by
the squares of the elements of $\pi_1(O_2)$. This is a finite index subgroup
of $\pi_1(O_2)$. Let $Z$ be the subgraph of ${\rm sing}(O_2)$ consisting of
the closure of the arcs and circles whose meridians
lie in $\langle \! \langle \mu \rangle \! \rangle \pi_1(O_2)^2$. 
Thus, $Z$ includes ${\rm sing}_p(O_2)$ and $V$. Let 
$O_3$ be the orbifold whose underlying manifold is the same
as that of $O_2$, and whose (possibly empty) singular locus is
${\rm sing}(O_2) - Z$. Note that ${\rm sing}(O_3)$
is a (possibly empty) collection of simple closed curves with order $2$. There is a natural
map $O_2 \rightarrow O_3$, which induces a homomorphism
$\pi_1(O_2) \rightarrow \pi_1(O_3)$. From the meridional presentations
for $\pi_1(O_2)$ and $\pi_1(O_3)$, we see that this homomorphism is surjective
and that its kernel
is normally generated by the meridians in $Z$, each of which lies in
$\langle \! \langle \mu \rangle \! \rangle \pi_1(O_2)^2$. Quotienting $\pi_1(O_2)$ by elements in $\pi_1(O_2)^2$
does not affect $d_2(O_2)$, and hence $d_2(O_3) \geq d_2(O_2) -1 \geq 9$.
By Lemma 2.1, $\pi_1(O_3)$ has a finite presentation $\langle X|R \rangle$ 
where $|R| - |X|$ is at most the number of components of ${\rm sing}(O_3)$,
which, by Lemma 2.2, is at most $d_2(O_3)$.
So,
$$d_2(O_3)^2/4 - |R| + |X| - d_2(O_3) \geq (d_2(O_3))^2/4 - 
2d_2(O_3) > 0,$$
because $d_2(O_3) \geq 9$. Hence, $\pi_1(O_3)$ is infinite,
by the Golod-Shafarevich inequality (see [16]). Therefore, so also is 
$\pi_1(O_3)^2$, since it is a finite index subgroup of $\pi_1(O_3)$.

Let $O_4$ be the covering space of $O_3$ corresponding to $\pi_1(O_3)^2$.
We claim that $O_4$ is a manifold. Consider any meridian of ${\rm sing}(O_3)$.
Its inverse image in $O_2$ is a meridian of ${\rm sing}(O_2) - Z$. 
Hence, it maps to a non-trivial element of $\pi_1(O_2)/
\langle \! \langle \mu \rangle \! \rangle \pi_1(O_2)^2$. Using meridional
presentations for $\pi_1(O_2)$ and $\pi_1(O_3)$, it is clear that
$\pi_1(O_2)/ \langle \! \langle \mu \rangle \! \rangle \pi_1(O_2)^2$
is isomorphic to $\pi_1(O_3)/\pi_1(O_3)^2$. Thus, no meridian of ${\rm sing}(O_3)$
lifts to a meridian of ${\rm sing}(O_4)$. Hence, ${\rm sing}(O_4)$ is empty.

Let $\tilde O$ be the cover of $O_2$
corresponding to $\langle \! \langle \mu \rangle \! \rangle \pi_1(O_2)^2$.
Clearly, $|\tilde O|$ is the same manifold as $|O_4|$. But we have
already shown that $O_4$ is a manifold. Hence, $\pi_1(|\tilde O|)
= \pi_1(O_4) = \pi_1(O_3)^2$, which we have established is infinite.
This is conclusion (3) of the theorem.

Since $\tilde O$ covers $O_2$, every arc and circle of ${\rm sing}(\tilde O)$
has order $2$ or $p$, which is (1) of the theorem. Moreover,
every vertex of ${\rm sing}(\tilde O)$ has local group isomorphic to $D$.
Note that the inverse image of $\overline{{\rm sing}_p(O_2)}$ in $\tilde O$ is
$\overline{{\rm sing}_p(\tilde O)}$, and that the inverse image of
$V$ lies in $\overline{{\rm sing}_2(\tilde O)}$. Indeed, the singular
locus of $\tilde O$ is the pre-image of $Z$. Hence, $\tilde O$
contains at least one vertex, giving (2) of the theorem.

Let $O_5$ (respectively, $O_6$) be the covering space of $O_1$ corresponding to
${\rm ker}(\phi)$ (respectively, $\phi^{-1}(C_1 \times \dots \times C_m)$). 
Then $O_5$ is a manifold, and $O_6$ is the quotient of $O_5$ by the
group of covering transformations $C_1 \times \dots \times C_m$, which is a product of copies
of ${\Bbb Z}/p$, and therefore ${\rm sing}(O_6)$ is a collection of
simple closed curves with order $p$. Note that $O_6$ is a regular
cover of $O_2$ with covering group $({\Bbb Z}/2)^m$.
Let $\tilde O'$ be the cover of $O_2$ corresponding to
$\pi_1(\tilde O) \cap \pi_1(O_6)$. This is a regular cover
of $\tilde O$ with degree a power of $2$. The singular locus is 
non-empty because ${\rm sing}_p(\tilde O)$ was non-empty. $\square$

The proof of Theorem 1.1 now proceeds almost as in the ${\Bbb Z}/2 \times {\Bbb Z}/2$
case. Since $|\tilde O|$ has infinite fundamental group, geometrisation
implies that either it has a finite cover with positive first Betti number
or it is hyperbolic. In the former case, $\pi_1(O)$ contains a
surface subgroup. So, we now suppose that $|\tilde O|$ is hyperbolic. Thus, by Theorem 1.3,
it has a sequence of infinite covers $|O_i|$ such that $h(|O_i|) \rightarrow 0$.
Let $O_i'$ be the cover of $\tilde O$ corresponding to $\pi_1(O_i)
\cap \pi_1(\tilde O')$. This is a regular cover of $O_i$ with degree a power of $2$,
at most that of the degree of $\tilde O' \rightarrow \tilde O$.

As in Section 3, pick a 1-vertex triangulation of $|\tilde O|$,
where the vertex is a singular point with local group $D$. Let $X_i$
be the inverse image of its 1-skeleton in $|O_i|$.
Then, we may find finite subsets $A_i$ of $V(X_i)$ such that
$|\partial A_i|/|A_i| \rightarrow 0$. These then
specify a closed orientable surface $|S_i|$ in $|O_i|$ bounding a
compact submanifold $|N_i|$. Let $S_i$ and $N_i$ be the underlying
orbifolds. Let $S'_i$ and $N'_i$ be their inverse images in $O_i'$.

Claim 1 in Section 3 gives that there are constants $k_0, k_1, k_2 > 0$, independent of $i$,
such that the number of 0-cells, 1-cells and 2-cells of $|S_i|$ is
at most $k_0 |\partial A_i|$, $k_1 |\partial A_i|$ and $k_2 |\partial A_i|$.
Also, the argument of Claim 2 gives there is a constant $k_3 > 0$, independent of
$i$, such that $d_p(S'_i) \leq k_3 |\partial A_i|$.
Let $k_4$ be the maximal number of intersection points between
any triangle or quadrilateral in $|S_i|$ and ${\rm sing}_p(O_i)$.
Now, there is a vertex of $N_i$ at each point of $A_i$.
Emanating from each of these vertices is an arc of $\overline{{\rm sing}_p(O_i)}$.
At most $k_4 k_2 |\partial A_i|$ of these arcs end on $S_i$.
The remainder lie entirely in $N_i$ and end at a vertex not in $A_i$. The inverse image of
any such arc in $N'_i$ is a collection of simple closed curves in
${\rm sing}_p(N'_i)$. Thus, we deduce that there are at least
$|A_i| - k_4 k_2 |\partial A_i|$ closed curves in ${\rm sing}_p(N'_i)$.
By Lemma 2.2, this is a lower bound for $d_p(N'_i)$.
So, if $|\partial A_i|/|A_i|$ is sufficiently small,
$d_p(N'_i) > d_p(S'_i)$. We now argue as in Section 3 to find
a non-trivial element of $H^1(O'_i; {\Bbb Z}/p)$ with trivial
image in $H^1(O'_i - N'_i; {\Bbb Z}/p)$. The $p$-fold
cover $O''_i$ of $O'_i$ dual to this class has at least $p$ ends. 

By Selberg's lemma, $O$ has a finite manifold cover $M$.
Let $M_i$ be the cover of $O$ corresponding to $\pi_1(M) \cap \pi_1(O''_i)$.
This is a finite cover of $O''_i$ and hence it also has at least $p$ ends.
It is a hyperbolic manifold, and hence irreducible.
So, by Lemma 3.3, its fundamental group contains a surface subgroup.
This is a surface subgroup of $\pi_1(O)$, thereby proving
Theorem 1.1. $\square$

\vfill\eject
\centerline{\caps 5. Comments on a theorem of Bowen}
\vskip 6pt

Crucial in the proof of Theorem 1.3 was Theorem 1.4, due to Bowen [2]. In this section,
we describe this result and provide an alternative interpretation of part of his proof.

For a Kleinian group $\G$, let
$D(\G)$ denote the Hausdorff dimension of its limit set. Bowen's result
(Theorem 1.4 and Remark 1 of [2]) is as follows.

\noindent {\bf Theorem 5.1.} {\sl Let $\G$ be a co-compact lattice in ${\rm PSL}(2, {\Bbb C})$.
Then $\G$ contains a sequence of finitely generated, free, convex co-compact
subgroups $\G_i$ such that $D(\G_i) \rightarrow 2$.}

In fact, Bowen proved rather more than this. He showed that $D(\G_i)$ can be
arranged to tend to any real number in the interval $[0,2]$. He also
dealt with lattices in ${\rm SO}(n,1)$, for any integer $n \geq 2$.
However, Theorem 5.1 is all that is needed to prove Theorem 1.3.

Theorem 5.1 is proved using $\epsilon$-perturbations, which are defined as follows.
Let $d$ be the left-invariant metric on ${\rm PSL}(2, {\Bbb C})$.
Let $F$ be a finitely generated free group, with a free generating
set $S$. Let $\phi \colon F \rightarrow {\rm PSL}(2, {\Bbb C})$ be a
homomorphism. Then, for $\epsilon > 0$, a map $\phi_\epsilon \colon F \rightarrow
{\rm PSL}(2, {\Bbb C})$ is an {\sl $\epsilon$-perturbation} of $\phi$ if
$d(\phi_\epsilon(fs), \phi_\epsilon(f) \phi(s)) \leq \epsilon$
for all $f \in F$ and $s \in S \cup S^{-1}$. Note that $\phi_\epsilon$ is
not assumed to be a homomorphism.

Let $\G$ be a lattice in ${\rm PSL}(2, {\Bbb C})$. Then 
$\phi_\epsilon \colon F \rightarrow {\rm PSL}(2, {\Bbb C})$ is {\sl virtually a homomorphism
into $\G$} if there exists a finite index subgroup $F'$ of $F$
such that $\phi_\epsilon(f' f) = \phi_\epsilon(f') \phi_\epsilon(f)$ for all $f' \in F'$
and $f \in F$, and where $\phi_\epsilon(F')$ is a subgroup of $\Gamma$.

Theorem 5.1 is a consequence of the following two results of Bowen, which
appear as Theorems 1.1 and 1.3 in [2].

\noindent {\bf Theorem 5.2.} {\sl Let $d$, $F$, $S$ and $\phi$ be as
above. Let $\G$ be a co-compact lattice in ${\rm PSL}(2, {\Bbb C})$.
Then, for any $\epsilon > 0$, there exists an $\epsilon$-perturbation
$\phi_\epsilon$ of $\phi$ that is virtually a homomorphism into
$\G$.}

\noindent {\bf Theorem 5.3.} {\sl Let $d$, $F$, $S$, $\phi$ and $\G$ be as
above. Suppose that $\phi$ is an injection onto a finitely generated,
convex co-compact subgroup of ${\rm PSL}(2, {\Bbb C})$.  For every $\epsilon > 0$, let $\phi_\epsilon$
be an $\epsilon$-perturbation of $\phi$ that is virtually
a homomorphism into $\G$. Let $F_\epsilon$ be any finite index
subgroup of $F$ such that $\phi_\epsilon(F_\epsilon)$ is a subgroup of $\G$
and, for all $f' \in F_\epsilon$ and all $f \in F$, $\phi_\epsilon(f'f)
= \phi_\epsilon(f') \phi_\epsilon(f)$. Then, for all sufficiently small $\epsilon$,
$\phi_\epsilon|F_\epsilon$ is an injection onto a convex co-compact subgroup.
Moreover, $D(\phi_\epsilon(F_\epsilon)) \rightarrow D(\phi(F))$ as $\epsilon \rightarrow 0$.}

Theorem 5.1 is a consequence of these results, via the following
argument. It is well known that ${\rm PSL}(2, {\Bbb C})$
contains finitely generated, free, convex co-compact subgroups $F$
with $D(F)$ arbitrarily close to $2$. (See Corollary 7.8 in [18] for example.)
Let $\phi \colon F \rightarrow {\rm PSL}(2, {\Bbb C})$
be the inclusion homomorphism. For each $\epsilon > 0$,
Theorem 5.2 gives an $\epsilon$-perturbation $\phi_\epsilon$ of $\phi$
that is a virtual homomorphism into $\G$. Let $F_\epsilon$ be the finite
index subgroup of $F$ in the definition of virtual homomorphism.
By Theorem 5.3, $\phi_\epsilon(F_\epsilon)$ is free, finitely
generated and convex co-compact, and $D(\phi_\epsilon(F_\epsilon)) \rightarrow D(F)$,
which was arbitrarily close to $2$.

The proof of Theorem 5.3 follows reasonably standard lines,
and we do not have anything to add to it. It is Theorem 5.2 that
has a more unusual proof. Bowen's argument uses symbolic dynamics
over the free group. Here we give an alternative proof, with essentially
the same mathematical content, but which is more directly geometric.

Let $d$, $F$, $S$, $\phi$ and $\G$ be as above. By passing to a finite index
subgroup if necessary, we may assume that $\G$ is torsion-free.
Hence, $M = \G \backslash {\Bbb H}^3$ is a manifold.
Let $\epsilon > 0$ be a given real number. Our aim is to find 
an $\epsilon$-perturbation of $\phi$. We may reduce $\epsilon$,
so that for any $g \in \G$ other than the identity, $d(g, {\rm id}) > \epsilon$. 
Let $\delta > 0$ be such that for all $s \in S$ and all $g_1, g_2$ in ${\rm PSL}(2,{\Bbb C})$
satisfying $d(g_1, {\rm id}) \leq \delta$ and $d(g_2, {\rm id}) \leq \delta$,
the inequalities $d(g_1 \phi(s) g_2, \phi(s)) \leq \epsilon$ and
$d(g_1 \phi(s^{-1}) g_2, \phi(s^{-1})) \leq \epsilon$ hold.

The aim is to construct
a finite graph $X$ that will specify the required $\epsilon$-perturbation of $\phi$.
In particular, the fundamental group of $X$ will be the finite index
subgroup $F'$ of $F$ in the definition of a virtual homomorphism.

Let $B = \G \backslash {\rm PSL}(2, {\Bbb C})$, the set of right cosets
of $\G$. This is the frame bundle over $M$. The metric $d$ descends to a metric
on $B$, which we will also call $d$. We partition $B$ into a finite
collection of sets $B_1, \dots, B_n$, each with diameter at most
$\delta$. Such a partition exists because $\G$ is co-compact.
We may assume that each $B_i$ is measurable with respect to Haar measure on $B$.
We may also assume that the boundary of each $B_i$ has measure zero.
In the interior of each $B_i$, we pick an element $\beta_i$.
This is a point in $M$, together with a frame at that point. We may
assume that $\beta_1$ is the identity coset in $B$. 

We now define a finite oriented graph $Y$.
There is a vertex of $Y$ for each $\beta_i$. Each edge of $Y$ has a unique label $s \in S$.
Two vertices $\beta_i$ and $\beta_j$ are joined by an $s$-labelled
edge, oriented from $\beta_i$ to $\beta_j$, if and only if there
is a $\beta'_i \in {\rm int} (B_i)$ and $\beta'_j \in {\rm int}(B_j)$ such that
$\beta_i' \phi(s) = \beta_j'$. Here, we are acting by right multiplication 
in ${\rm PSL}(2, {\Bbb C})$. Thus, the frame $\beta_j'$ is obtained from the
frame $\beta_i'$ by post-composing with the isometry $\phi(s)$.

For each $s$-labelled edge $e$ in $Y$, we now define $\psi(e)$ in
${\rm PSL}(2, {\Bbb C})$, which will be an approximation to $\phi(s)$.
Let $\beta_i$ and $\beta_j$ be the initial and terminal vertices
of $e$. By definition of the edge $e$, there is some $\beta'_i \in {\rm int}(B_i)$
and $\beta'_j \in {\rm int}(B_j)$ such that $\beta'_i \phi(s) = \beta'_j$.
Since the sets $B_i$ and $B_j$ each have diameter at most $\delta$,
there are elements $g_i$ and $g_j$ in ${\rm PSL}(2, {\Bbb C})$
such that $\beta_i g_i = \beta'_i$, $\beta_j g_j = \beta'_j$,
$d(g_i, {\rm id}) \leq \delta$ and $d(g_j, {\rm id}) \leq \delta$.
Define $\psi(e)$ to be $g_i \phi(s) g_j^{-1}$. Then $\beta_i \psi(e)
= \beta_j$. Also, by definition of $\delta$, $d(\psi(e), \phi(s)) \leq \epsilon$
and $d(\psi(e)^{-1}, \phi(s^{-1})) \leq \epsilon$.

One can view each vertex $\beta_i$ of $Y$ as placed at the corresponding
point of $B$ (or $M$), and each edge $e$ running in the direction of $\psi(e)$.
(See Figure 2.)

\vskip 18pt
\centerline{
\epsfxsize=3in
\epsfbox{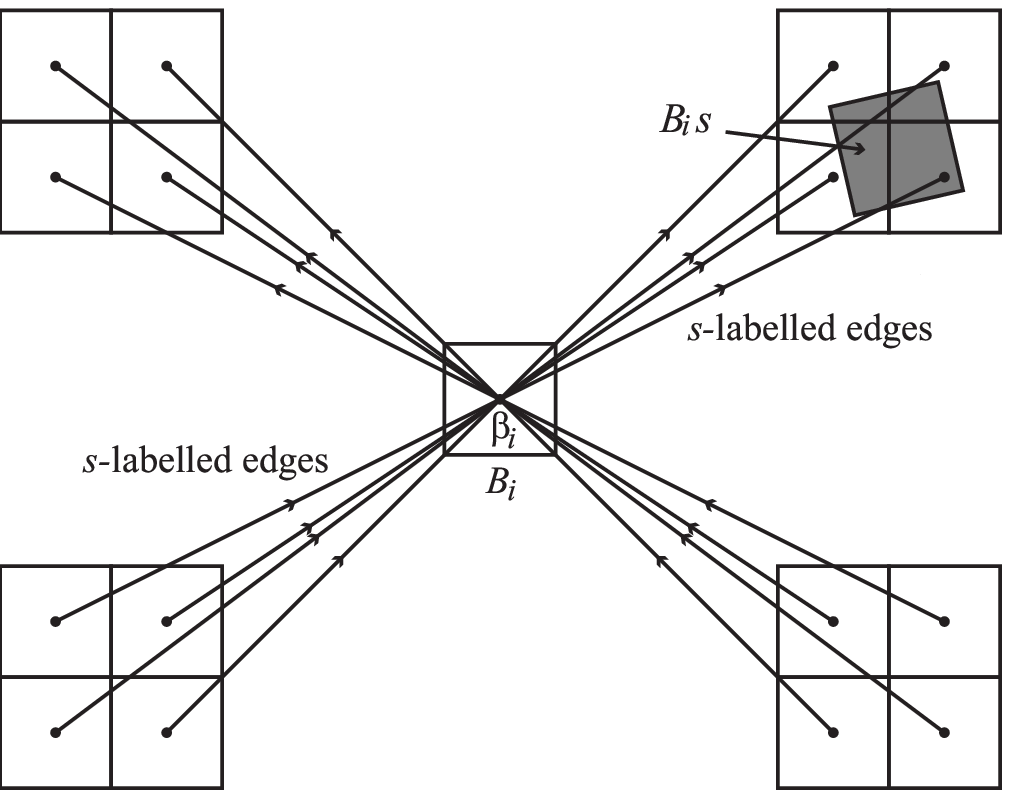}
}
\vskip 6pt
\centerline{Figure 2: The link of a vertex $\beta_i$}

Define a {\sl weighting} on $Y$ to be functions $w \colon V(Y) \rightarrow (0,\infty)$
and $w \colon E(Y) \rightarrow (0, \infty)$, where for each $v \in V(Y)$ and
$e \in E(Y)$, $w(v)$ and $w(e)$ are known as the {\sl weights} of $v$ and $e$.
We require that, for each $s \in S$ and $v \in V(Y)$,
the sum of the weights of the $s$-labelled edges emanating from $v$ is equal to
the weight of $v$, and that the sum of the weights of the $s$-labelled edges
entering $v$ is also the weight of $v$. Note that the conditions specifying a weighting
on $Y$ are a finite collection of linear equations with integer coefficients.

We claim that there is a weighting $w$ on $Y$. Let $\mu$ be Haar measure on
$B$. For each vertex $\beta_i$ of
$Y$, define $w(\beta_i)$ to be $\mu(B_i)$. Each edge $e$ of $Y$ has a 
label $s \in S$. Let $\beta_i$ and $\beta_j$ be its initial and
terminal vertices. 
Set $w(e)$ to be 
$$\mu(\{ \beta'_i \in B_i : \beta_i' \phi(s) \in B_j \}).$$
Note that $w(e)$ is positive,
because, by the definition of the edges of $Y$, the set
$\{ \beta'_i \in B_i : \beta_i' \phi(s) \in B_j \}$ has non-empty interior.
Consider the set of $s$-labelled
edges emanating from $\beta_i$. The sum of the weights of these
edges is clearly $\mu(B_i)$, which is $w(\beta_i)$.
Now consider the sum of the weights of the $s$-labelled edges
entering $\beta_i$. Each such edge starts at some $\beta_j$.
It has weight
$$\eqalign{
\mu(\{ \beta'_j \in B_j : \beta_j' \phi(s) \in B_i \}) &= 
\mu (\{ \beta'_j \phi(s) : \beta'_j \in B_j \hbox { and } \beta_j' \phi(s) \in B_i \}) \cr
&= \mu ( \{ \beta'_i \in B_i : \beta_i' \phi(s^{-1}) \in B_j \}).}$$
The first equality is a consequence of the invariance of Haar measure
under right multiplication.
Clearly, the sum of the weights over all $B_j$ is again
$\mu(B_i)$, which is $w(\beta_i)$. Thus,
$w$ is indeed a weighting as claimed.

Since the defining conditions for a weighting are a finite collection of linear equations with
integer coefficients, the existence of a weighting implies the
existence of a weighting $w'$ with integer weights.

We use $w'$ to construct the finite graph $X$. Associated with each
vertex $v$ of $Y$, there will be $w'(v)$ vertices of $X$. Similarly,
for each edge $e$ of $Y$, there will be $w'(e)$ edges of $X$.
Formally, this defines functions $h \colon V(X) \rightarrow V(Y)$
and $h \colon E(X) \rightarrow E(Y)$ such that, for every $v \in V(Y)$,
$|h^{-1}(v)| = w'(v)$ and, for every $e \in E(Y)$, $|h^{-1}(e)| = w'(e)$.
We need to specify the initial
and terminal vertices of edges of $X$. Consider a vertex $v \in V(Y)$. 
Its inverse image in $V(X)$ is $w'(v)$
vertices. For any $s \in S$, the total weight of the $s$-labelled edges emanating
from $v$ is also $w'(v)$. Thus, their inverse images in $Y$ form a
set of $w'(v)$ edges. Pick an arbitrary bijection from this set
of edges to $h^{-1}(v)$. Declare that these are the initial
vertices of these edges. Similarly, the total weight of the $s$-labelled
edges entering $v$ is also $w'(v)$. Their inverse images in $X$ again form
a set of $w'(v)$ edges. Pick an arbitrary bijection from this
set to $h^{-1}(v)$, and declare that these are the terminal
vertices of these edges. Do this for every $v \in V(Y)$ and $s \in S$
to create a graph. Note that this graph need not be connected. If it is not, we discard
all but one component. We may arrange that some vertex $b$ in this
component maps to $\beta_1$.
Thus, we have constructed the graph $X$ together with a continuous mapping
$h \colon X \rightarrow Y$ that sends vertices to vertices, and
sends the interiors of edges homeomorphically to the interiors of edges.

The edges of $X$ are oriented and labelled by the elements of $S$.
Each vertex has exactly one $s$-labelled edge pointing into it,
and one $s$-labelled edge pointing out, for each $s \in S$. Thus,
$X$ is naturally a covering space of the wedge of $|S|$ circles, with labels
in $S$. This wedge of circles has fundamental group $F$,
and so the fundamental group $\pi_1(X,b)$ corresponds to a subgroup
$F'$ of $F$. This has finite index in $F$ because $X$ is a finite graph.

We can now define the $\epsilon$-approximation $\phi_\epsilon$.
Let $f$ be any element of $F$, and let $s_1^{k_1} \dots s_n^{k_n}$ (where each $k_j \in \{ -1,1 \}$) 
be any word in the generators representing $f$. This specifies
a path $e_1^{k_1} \dots e_n^{k_n}$ in $X$ starting at $b$, where $e_i^{-1}$ denotes
the reverse of $e_i$. Define $\phi_\epsilon(f)$
to be $\psi(h(e_1))^{k_1} \dots \psi(h(e_n))^{k_n}$. We need to check that this
is independent of the choice of word representing $f$. But any
two such words differ by a sequence of moves, each of which
inserts or removes $ss^{-1}$, for some $s \in S \cup S^{-1}$. This
changes the path $e_1^{k_1} \dots e_n^{k_n}$ by the insertion or removal
of $ee^{-1}$ for some edge $e$. This leaves $\phi_\epsilon(f)$ unchanged
and therefore well-defined.

We now check that $\phi_\epsilon$ is an $\epsilon$-approximation. Consider
$f \in F$ and $s \in S$. Let $e_1^{k_1} \dots e_n^{k_n} e^{k}$ be a path in $X$ starting at $b$
representing $fs^k$, where $k \in \{-1,1 \}$. Then, $\phi_\epsilon(fs^k)  = \phi_\epsilon(f) \psi(h(e))^k$,
and so $d(\phi_\epsilon(fs^k), \phi_\epsilon(f) \phi(s^k)) = d(\psi(h(e))^k, \phi(s^k)) \leq \epsilon$.
Thus, $\phi_\epsilon$ is indeed an $\epsilon$-approximation.

We also claim that $\phi_\epsilon$ is virtually a homomorphism into $\Gamma$.
Recall that $F'$ is those elements of $F$ that give rise to loops in $X$ based at $b$,
which is a finite index subgroup of $F$.
Then, by construction, for all $f' \in F'$, $\beta_1 \phi_\epsilon(f') = \beta_1$.
In other words, $\phi_\epsilon(f') \in \G$. Suppose that $f' \in F'$
and $f \in F$. Let $(e'_1)^{k_1'} \dots (e'_m)^{k'_m}$ be a loop based at $b$ representing $f'$,
and let $e_1^{k_1} \dots e_n^{k_n}$ be a path starting at $b$ representing $f$. Then
$(e'_1)^{k'_1} \dots (e'_m)^{k'_m} e_1^{k_1} \dots e_n^{k_n}$ is a path starting at $b$ representing
$f'f$, because $(e'_1)^{k_1'} \dots (e'_m)^{k'_m}$ is a loop. Hence,
$$\phi_\epsilon(f'f) = \psi(h(e'_1))^{k'_1} \dots \psi(h(e'_m))^{k'_m} \psi(h(e_1))^{k_1} \dots \psi(h(e_n))^{k_n}
= \phi_\epsilon(f') \phi_\epsilon(f),$$
as required.

This completes the proof of Theorem 5.2. $\square$

\vskip 18pt
\centerline{\caps References}
\vskip 6pt

\item{1.} {\caps M. Bestvina}, {\sl Questions in geometric group theory},
available at \hfill\break http://www.math.utah.edu/$\sim$bestvina
\item{2.} {\caps L. Bowen}, {\sl Free groups in lattices}, Preprint, available
at arxiv:0802.0185
\item{3.} {\caps R. Brooks}, {\sl The fundamental group and the spectrum of the
Laplacian}, Comment. Math. Helv. 56 (1981) 581--598.
\item{4.} {\caps J. Cheeger}, {\sl A lower bound for the smallest eigenvalue of
the Laplacian}, Problems in analysis (Papers dedicated to Salomon Bochner, 1969),
pp 195--199, Princeton University Press.
\item{5.} {\caps L. Clozel}, {\sl On the cuspidal cohomology of arithmetic subgroups
of ${\rm SL}(2n)$ and the first Betti number of arithmetic 3-manifolds},
Duke Math. J. 55 (1987) 475--486.
\item{6.} {\caps D. Cooper, C. Hodgson, S. Kerckhoff,}
{\sl Three-dimensional orbifolds and cone-manifolds.} MSJ Memoirs, 5. 
Mathematical Society of Japan, Tokyo, 2000.
\item{7.} {\caps D. Cooper, D. Long, A. Reid,} 
{\sl Essential closed surfaces in bounded $3$-manifolds.}
J. Amer. Math. Soc. 10 (1997) 553--563.
\item{8.} {\caps N. Dunfield, W. Thurston}, {\sl The virtual Haken conjecture:
experiments and examples}, Geom. Topol. 7 (2003) 399-441.
\item{9.} {\caps J. Hempel}, {\sl $3$-manifolds}, Ann. of Math. Studies,
Princeton (1976).
\item{10.} {\caps J-P. Labesse, J. Schwermer}, {\sl On liftings and cusp cohomology
of arithmetic groups}, Invent. Math. 83 (1986) 383--401.
\item{11.} {\caps M. Lackenby}, {\sl Heegaard splittings,
the virtually Haken conjecture and Property ($\tau$)},
Invent. Math. 164 (2006) 317--359.
\item{12.} {\caps M. Lackenby}, {\sl Large groups, Property $(\tau)$ and
the homology growth of subgroups}, Math. Proc. Camb. Phil. Soc. 146 (2009) 625--648.
\item{13.} {\caps M. Lackenby}, {\sl Covering spaces of 3-orbifolds},
Duke Math J. 136 (2007) 181--203.
\item{14.} {\caps M. Lackenby, D. Long, A. Reid}, 
{\sl Covering spaces of arithmetic 3-orbifolds}, Int. Math. Res. Not. IMRN,
(2008).
\item{15.} {\caps M. Lackenby, D. Long, A. Reid}, {\sl
LERF and the Lubotzky-Sarnak conjecture}, Geom. Topol. 12 (2008) 2047--2056. 
\item{16.} {\caps A. Lubotzky}, {\sl Group presentations, $p$-adic analytic
groups and lattices in ${\rm SL}_2({\Bbb C})$}, Ann. Math. 118 (1983) 115--130.
\item{17.} {\caps A. Lubotzky}, {\sl Eigenvalues of the 
Laplacian, the first Betti number and the congruence subgroup
problem,} Ann. Math. 144 (1996) 441--452.
\item{18.} {\caps C. McMullen, } {\sl Hausdorff dimension and conformal dynamics. 
I. Strong convergence of Kleinian groups.} J. Differential Geom. 51 (1999) 471--515.
\item{19.} {\caps G. Perelman,} {\sl The entropy formula for the 
Ricci flow and its geometric applications,} Preprint,
available at arxiv:math.DG/0211159
\item{20.} {\caps G. Perelman,} {\sl  Ricci flow with surgery on three-manifolds,}
Preprint, available at arxiv:math.DG/0303109
\item{21.} {\caps G. Perelman,} {\sl Finite extinction time for the 
solutions to the Ricci flow on certain three-manifolds,} Preprint,
available at arxiv:math.DG/0307245 
\item{22.} {\caps P. Scott}, {\sl Compact submanifolds of $3$-manifolds},
J. London Math. Soc. 7 (1973) 246--250.
\item{23.} {\caps P. Scott}, {\sl The geometries of $3$-manifolds},
Bull. London Math. Soc. 15 (1983) 401--487.
\item{24.} {\caps D. Sullivan,} {\sl Entropy, Hausdorff measures old and new,
and limit sets of geometrically finite Kleinian groups,} Acta. Math. 153
(1984) 259--277.
\item{25.} {\caps W. Thurston}, {\sl Three-dimensional manifolds, Kleinian groups
and hyperbolic geometry}, Bull. Amer. Math. Soc. 6 (1982) 357--381.

\vskip 12pt

\+ Mathematical Institute, University of Oxford, \cr
\+ 24-29 St Giles', Oxford OX1 3LB, United Kingdom. \cr

\end